\documentclass[12pt]{amsart}
\usepackage{amscd,amssymb,amsthm}
\usepackage{graphicx}

\newtheorem{theorem}{Theorem}[section]
\newtheorem{remark}{Remark}[section]
\newtheorem{example}{Example}[section]

\newtheorem{proposition}{Proposition}[section]
\numberwithin{equation}{section}

\title[Refinement Inequalities]{Refinement Inequalities Among Symmetric
Divergence Measures}

\author{Inder Jeet Taneja}
\address{Inder Jeet Taneja \\
Departamento de Matem\'{a}tica\\
Universidade Federal
de Santa Catarina\\
88.040-900 Florian\'{o}polis, SC, Brazil}

\email{taneja@mtm.ufsc.br}

\urladdr{http://www.mtm.ufsc.br/$\sim $taneja}


\begin{document}
\begin{abstract}
There are three classical divergence measures in the literature on
information theory and statistics, namely,
Jeffryes-Kullback-Leiber's \textit{J-divergence},
Sibson-Burbea-Rao's \textit{Jensen-Shannon divegernce} and Taneja's
\textit{arithemtic - geometric mean divergence}. These bear an
interesting relationship among each other and are based on
logarithmic expressions. The divergence measures like
\textit{Hellinger discrimination}, \textit{symmetric }$\chi ^2 -
$\textit{divergence}, and \textit{triangular discrimination} are not
based on logarithmic expressions. These six divergence measures are
symmetric with respect to probability distributions. In this paper
some interesting inequalities among these symmetric divergence
measures are studied. Refinement over these inequalities is also
given. Some inequalities due to Dragomir et al. \cite{dsb} are also
improved.
\end{abstract}

\maketitle

\section{Introduction}
Let
\[
\Gamma _n = \left\{ {P = (p_1 ,p_2 ,...,p_n )\left| {p_i >
0,\sum\limits_{i = 1}^n {p_i = 1} } \right.} \right\}, \, n
\geqslant 2,
\]

\noindent be the set of all complete finite discrete probability
distributions. For all $P,Q \in \Gamma _n $, the following
measures are well known in the literature on information theory
and statistics:

\bigskip
\textbf{$\bullet$ Hellinger Discrimination}
\begin{equation}
\label{eq1} h(P\vert \vert Q) = 1 - B(P\vert \vert Q) =
\frac{1}{2}\sum\limits_{i = 1}^n {(\sqrt {p_i } - \sqrt {q_i }
)^2},
\end{equation}

\noindent where
\begin{equation}
\label{eq2}
B(P\vert \vert Q) = \sqrt {p_i q_i } ,
\end{equation}

\noindent is the well-known Bhattacharyya \cite{bha}
\textit{coefficient}.

\bigskip
\textbf{$\bullet$ Triangular Discrimination}
\begin{equation}
\label{eq3}
\Delta (P\vert \vert Q) = 2\left[ {1 - W(P\vert \vert Q)} \right] =
\sum\limits_{i = 1}^n {\frac{(p_i - q_i )^2}{p_i + q_i }} ,
\end{equation}

\noindent where
\begin{equation}
\label{eq4}
W(P\vert \vert Q) = \sum\limits_{i = 1}^n {\frac{2p_i q_i }{p_i + q_i }} ,
\end{equation}

\noindent is the well-known \textit{harmonic mean divergence}.

\bigskip
\textbf{$\bullet$ Symmetric Chi-square Divergence}
\begin{equation}
\label{eq5}
\Psi (P\vert \vert Q) = \chi ^2(P\vert \vert Q) + \chi ^2(Q\vert \vert P) =
\sum\limits_{i = 1}^n {\frac{(p_i - q_i )^2(p_i + q_i )}{p_i q_i }} ,
\end{equation}

\noindent where
\begin{equation}
\label{eq6}
\chi ^2(P\vert \vert Q) = \sum\limits_{i = 1}^n {\frac{(p_i - q_i )^2}{q_i
}} = \sum\limits_{i = 1}^n {\frac{p_i^2 }{q_i } - 1} ,
\end{equation}

\noindent is the well-known $\chi ^2 - $\textit{divergence}
(Pearson \cite{pea})

\bigskip
\textbf{$\bullet$ J-Divergence}
\begin{equation}
\label{eq7}
J(P\vert \vert Q) = \sum\limits_{i = 1}^n {(p_i - q_i )\ln (\frac{p_i }{q_i
})} .
\end{equation}

\bigskip
\textbf{$\bullet$ Jensen-Shannon Divergence}
\begin{equation}
\label{eq8}
I(P\vert \vert Q) = \frac{1}{2}\left[ {\sum\limits_{i = 1}^n {p_i \ln \left(
{\frac{2p_i }{p_i + q_i }} \right) + } \sum\limits_{i = 1}^n {q_i \ln \left(
{\frac{2q_i }{p_i + q_i }} \right)} } \right].
\end{equation}

\bigskip
\textbf{$\bullet$ Arithmetic-Geometric Mean Divergence}
\begin{equation}
\label{eq9}
T(P\vert \vert Q) = \sum\limits_{i = 1}^n {\left( {\frac{p_i + q_i }{2}}
\right)\ln \left( {\frac{p_i + q_i }{2\sqrt {p_i q_i } }} \right)} .
\end{equation}

After simplification, we can write
\begin{equation}
\label{eq10}
J(P\vert \vert Q) = 4\left[ {I(P\vert \vert Q) + T(P\vert \vert Q)}
\right].
\end{equation}

\bigskip
The measures $I(P\vert \vert Q)$, $J(P\vert \vert Q)$ and
$T(P\vert \vert Q)$ can be written as
\begin{align}
\label{eq11} J(P\vert \vert Q) & = K(P\vert \vert Q) + K(Q\vert
\vert P),\\
\label{eq12} I(P\vert \vert Q) & = \frac{1}{2}\left[ {K\left(
{P\vert \vert \frac{P + Q}{2}} \right) + K\left( {Q\vert \vert
\frac{P + Q}{2}} \right)} \right], \\
\intertext{and} \label{eq13}
T(P\vert \vert Q) & = \frac{1}{2}\left[ {K\left( {\frac{P +
Q}{2}\vert \vert P} \right) + K\left( {\frac{P + Q}{2}\vert \vert
Q} \right)} \right],
\end{align}

\noindent where
\begin{equation}
\label{eq14}
K(P\vert \vert Q) = \sum\limits_{i = 1}^n {p_i \log \left( {\frac{p_i }{q_i
}} \right)} ,
\end{equation}

\noindent is the well known Kullback-Leibler \cite{kul}
\textit{relative information}.

We call the measures given in (\ref{eq1}), (\ref{eq3}), (\ref{eq5}),
(\ref{eq7}), (\ref{eq9}) and (\ref{eq10}) as \textit{symmetric
divergence measures}, since they are symmetric with respect to the
probability distributions $P$ and $Q$. The measure (\ref{eq1}) is
due to Hellinger \cite{hel}. The measure (\ref{eq5}) is due to
Dragomir et al. \cite{dsb}, and recently has been studied by Taneja
\cite{tan4}. The measure (\ref{eq7}) is due to Jeffreys \cite{jef},
and later Kullback-Leibler \cite{kul} studied it extensively. Some
times it is called as Jeffreys-Kullback-Leibler's
\textit{J-divergence}. The measure (\ref{eq8}) is due to Sibson
\cite{sib}, and later Burbea and Rao \cite{bur1, bur2} studied it
extensively. Initially, it was called as \textit{information
radius}, but now a days it is famous as \textit{Jensen-Shannon
divegence}. The measure (\ref{eq9}) is due to Taneja \cite{tan4},
and is known by \textit{arithmetic-geometric mean divergence}. For
one parametric generalizations of the measures given above refer to
Taneja \cite{tan6, tan7}. A general study of information and
divergence measures and their generalizations can be seen in Taneja
\cite{tan1,tan2,tan3}.

In this paper our aim is to obtain an inequality its improvement in
terms of above symmetric divergence measures. This we shall do by
the application of some properties of Csisz\'{a}r's
$f-$\textit{divergence}.

\section{Csisz\'{a}r's $f-$Divergence}
Given a function $f:[0,\infty ) \to \mathbb{R}$, the
\textit{f-divergence} measure introduced by Csisz\'{a}r's
\cite{csi1} is given by
\begin{equation}
\label{eq19}
C_f (P\vert \vert Q) =
\sum\limits_{i = 1}^n {q_i f\left( {\frac{p_i }{q_i }} \right)} ,
\end{equation}

\noindent for all $P,Q \in \Gamma _n $.

The following theorem is well known in the literature.

\begin{theorem} \label{the21} (Csisz\'{a}r's \cite{csi1, csi2}). If the
function $f$ is convex and normalized, i.e., $f(1) = 0$, then the
$f-$divergence, $C_f (P\vert \vert Q)$ is nonnegative and convex in
the pair of probability distribution $(P,Q) \in \Gamma _n \times
\Gamma _n $.
\end{theorem}

Recently, Taneja \cite{tan5, tan7} established the following
property of the measure (\ref{eq19}).

\begin{theorem} \label{the22} (Taneja \cite{tan4}). Let $f_1 ,f_2 :I
\subset \mathbb{R}_ + \to \mathbb{R}$ two generating mappings are
normalized, i.e., $f_1 (1) = f_2 (1) = 0$ and satisfy the
assumptions:

(i) $f_1 $ and $f_2 $ are twice differentiable on $(a,b)$;

(ii) there exists the real constants $m,M$ such that $m < M$ and
\begin{equation}
\label{eq20} m \leqslant \frac{f_1 ^{\prime \prime }(x)}{f_2
^{\prime \prime }(x)} \leqslant M, \, f_2 ^{\prime \prime }(x) >
0, \, \forall x \in (a,b),
\end{equation}

\noindent then we have the inequalities:
\begin{equation}
\label{eq21}
m\mbox{ }C_{f_2 } (P\vert \vert Q) \leqslant C_{f_1 } (P\vert \vert Q)
\leqslant M\mbox{ }C_{f_2 } (P\vert \vert Q).
\end{equation}
\end{theorem}

\begin{proof} Let us consider the functions $\eta _{m.s} (
\cdot )$ and $\eta _{M.s} ( \cdot )$ given by
\begin{equation}
\label{eq22}
\eta _m (x) = f_1 (x) - m\mbox{ }f_2 (x),
\end{equation}

\noindent and
\begin{equation}
\label{eq23}
\eta _M (x) = M\mbox{ }f_2 (x) - f_1 (x),
\end{equation}

\noindent respectively, where $m$ and $M$ are as given by
(\ref{eq20}).

Since $f_1 (x)$ and $f_2 (x)$ are normalized, i.e., $f_1 (1) = f_2
(1) = 0$, then $\eta _m ( \cdot )$ and $\eta _M ( \cdot )$ are
also normalized, i.e., $\eta _m (1) = 0$ and $\eta _M (1) = 0$.
Also, the functions $f_1 (x)$ and $f_2 (x)$ are twice
differentiable. Then in view of (\ref{eq20}), we have
\begin{equation}
\label{eq24}
{\eta }''_m (x) = f_1 ^{\prime \prime }(x) - m\mbox{ }f_2 ^{\prime \prime
}(x)
 = f_2 ^{\prime \prime }(x)\left( {\frac{f_1 ^{\prime \prime }(x)}{f_2
^{\prime \prime }(x)} - m} \right) \geqslant 0,
\end{equation}

\noindent and
\begin{equation}
\label{eq25}
{\eta }''_M (x) = M\mbox{ }f_2 ^{\prime \prime }(x) - f_1 ^{\prime \prime
}(x)
 = f_2 ^{\prime \prime }(x)\left( {M - \frac{f_1 ^{\prime \prime }(x)}{f_2
^{\prime \prime }(x)}} \right) \geqslant 0,
\end{equation}

\noindent for all $x \in (r,R)$.

In view of (\ref{eq24}) and (\ref{eq25}), we can say that the
functions $\eta _m ( \cdot )$ and $\eta _M ( \cdot )$ given by
(\ref{eq22}) and (\ref{eq23}) respectively, are convex on $(r,R)$.

According to Theorem \ref{the21}, we have
\begin{equation}
\label{eq26}
C_{\eta _m } (P\vert \vert Q) = C_{f_1 - mf_2 } (P\vert \vert Q) = C_{f_1 }
(P\vert \vert Q) - m\mbox{ }C_{f_2 } (P\vert \vert Q) \geqslant 0,
\end{equation}

\noindent and
\begin{equation}
\label{eq27}
C_{\eta _M } (P\vert \vert Q) = C_{Mf_2 - f_1 } (P\vert \vert Q) = M\mbox{
}C_{f_2 } (P\vert \vert Q) - C_{f_1 } (P\vert \vert Q) \geqslant 0.
\end{equation}

Combining (\ref{eq26}) and (\ref{eq27}) we have the proof of
(\ref{eq21}).
\end{proof}

Now, based on Theorem \ref{the21}, we shall give below the
\textit{convexity} and \textit{nonnegativity} of the
\textit{symmetric divergence measures} given in Section 1.

\begin{example} \label{exa21} (\textit{Hellinger discrimination}). Let
us consider
\begin{equation}
\label{eq28} f_h (x) = \frac{1}{2}(\sqrt x - 1)^2, \, x \in
(0,\infty ),
\end{equation}

\noindent in (\ref{eq19}), then  we have $C_f (P\vert \vert Q) =
h(P\vert \vert Q),$ where $h(P\vert \vert Q)$ is as given by
(\ref{eq1}).

Moreover,
\[
{f}'_h (x) = \frac{\sqrt x - 1}{2\sqrt x },
\]

\noindent and
\begin{equation}
\label{eq29}
{f}''_h (x) = \frac{1}{4x\sqrt x }.
\end{equation}

Thus we have ${f}''_h (x) > 0$ for all $x > 0$, and hence, $f_h
(x)$ is \textit{strictly convex} for all $x > 0$. Also, we have
$f_h (1) = 0$. In view of this we can say that the
\textit{Hellinger discrimination }given by (\ref{eq1}) is
\textit{nonnegative} and \textit{convex} in the pair of
probability distributions $(P,Q) \in \Gamma _n \times \Gamma _n $.
\end{example}

\begin{example} \label{exa22} (\textit{Triangular discrimination}).
Let us consider
\begin{equation}
\label{eq30} f_\Delta (x) = \frac{(x - 1)^2}{x + 1}, \, x \in
(0,\infty ),
\end{equation}

\noindent in (\ref{eq19}), then  we have $C_f (P\vert \vert Q) =
\Delta (P\vert \vert Q),$ where $\Delta (P\vert \vert Q)$ is as
given by (\ref{eq3}).

Moreover,
\[
{f}'_\Delta (x) = \frac{(x - 1)(x + 3)}{(x + 1)^2},
\]

\noindent and
\begin{equation}
\label{eq31}
{f}''_\Delta (x) = \frac{8}{(x + 1)^3}.
\end{equation}

Thus we have ${f}''_\Delta (x) > 0$ for all $x > 0$, and hence,
$f_\Delta (x)$ is \textit{strictly convex} for all $x > 0$. Also,
we have $f_\Delta (1) = 0$. In view of this we can say that the
\textit{triangular discrimination }given by (\ref{eq3}) is
\textit{nonnegative} and \textit{convex} in the pair of
probability distributions $(P,Q) \in \Gamma _n \times \Gamma _n $.
\end{example}

\begin{example} \label{exa23} (\textit{Symmetric chi-square
divergence}). Let us consider
\begin{equation}
\label{eq32} f_\Psi (x) = \frac{(x - 1)^2(x + 1)}{x}, \, x \in
(0,\infty ),
\end{equation}

\noindent in (\ref{eq19}), then  we have $C_f (P\vert \vert Q) =
\Psi (P\vert \vert Q),$ where $\Psi (P\vert \vert Q)$ is as given
by (\ref{eq5}).

Moreover,
\[
{f}'_\Psi (x) = \frac{(x - 1)(2x^2 + x + 1)}{x^2},
\]

\noindent and
\begin{equation}
\label{eq33}
{f}''_\Psi (x) = \frac{2(x^3 + 1)}{x^3}.
\end{equation}

Thus we have ${f}''_\Psi (x) > 0$ for all $x > 0$, and hence,
$f_\Psi (x)$ is \textit{strictly convex} for all $x > 0$. Also, we
have $f_\Psi (1) = 0$. In view of this we can say that the
\textit{symmetric chi-square divergence} given by (\ref{eq5}) is
\textit{nonnegative} and \textit{convex} in the pair of
probability distributions $(P,Q) \in \Gamma _n \times \Gamma _n $.
\end{example}

\begin{example} \label{exa24} (\textit{J-divergence}). Let us consider
\begin{equation}
\label{eq34} f_J (x) = (x - 1)\ln x, \, x \in (0,\infty ),
\end{equation}

\noindent in (\ref{eq19}), then  we have $C_f (P\vert \vert Q) =
J(P\vert \vert Q),$ where $J(P\vert \vert Q)$ is as given by
(\ref{eq7}).

Moreover,
\[
{f}'_J (x) = 1 - x^{ - 1} + \ln x,
\]

\noindent and
\begin{equation}
\label{eq35}
{f}''_J (x) = \frac{x + 1}{x^2}.
\end{equation}

Thus we have ${f}''_J (x) > 0$ for all $x > 0$, and hence, $f_J
(x)$ is \textit{strictly convex} for all $x > 0$. Also, we have
$f_J (1) = 0$. In view of this we can say that the
\textit{J-divergence }given by (\ref{eq7}) is \textit{nonnegative}
and \textit{convex} in the pair of probability distributions
$(P,Q) \in \Gamma _n \times \Gamma _n $.
\end{example}

\begin{example} \label{exa25}(\textit{JS-divergence}). Let us
consider
\begin{equation}
\label{eq36} f_I (x) = \frac{x}{2}\ln x + \frac{x + 1}{2}\ln
\left( {\frac{2}{x + 1}} \right), \, x \in (0,\infty ),
\end{equation}

\noindent in (\ref{eq19}), then  we have $C_f (P\vert \vert Q) =
I(P\vert \vert Q),$ where $I(P\vert \vert Q)$ is as given by
(\ref{eq8}).

Moreover,
\[
{f}'_I (x) = \frac{1}{2}\ln \left( {\frac{2x}{x + 1}} \right),
\]

\noindent and
\begin{equation}
\label{eq37}
{f}''_I (x) = \frac{1}{2x(x + 1)}.
\end{equation}

Thus we have ${f}''_I (x) > 0$ for all $x > 0$, and hence, $f_I
(x)$ is \textit{strictly convex} for all $x > 0$. Also, we have
$f_I (1) = 0$. In view of this we can say that the
\textit{JS-divergence }given by (\ref{eq8}) is
\textit{nonnegative} and \textit{convex} in the pair of
probability distributions $(P,Q) \in \Gamma _n \times \Gamma _n $.
\end{example}

\begin{example} \label{exa26} (\textit{AG-Divergence}). Let us
consider
\begin{equation}
\label{eq38} f_T (x) = \left( {\frac{x + 1}{2}} \right)\ln \left(
{\frac{x + 1}{2\sqrt x }} \right), \, x \in (0,\infty ),
\end{equation}

\noindent in (\ref{eq19}), then we have $C_f (P\vert \vert Q) =
T(P\vert \vert Q),$ where $T(P\vert \vert Q)$ is as given by
(\ref{eq9}).

Moreover,
\[
{f}'_T (x) = \frac{1}{4}\left[ {1 - x^{ - 1} + 2\ln \left( {\frac{x +
1}{2\sqrt x }} \right)} \right],
\]

\noindent and
\begin{equation}
\label{eq39}
{f}''_T (x) = \frac{x^2 + 1}{4x^2(x + 1)}.
\end{equation}

Thus we have ${f}''_T (x) > 0$ for all $x > 0$, and hence, $f_T
(x)$ is \textit{strictly convex} for all $x > 0$. Also, we have
$f_T (1) = 0$. In view of this we can say that the
\textit{AG-divergence} given by (\ref{eq9}) is
\textit{nonnegative} and \textit{convex} in the pair of
probability distributions $(P,Q) \in \Gamma _n \times \Gamma _n $.
\end{example}

\section{Inequalities Among the Measures}

In this section we shall apply the Theorem \ref{the22} to obtain
inequalities among the measures given in Section 1. We have
considered only the symmetric measures given in (\ref{eq1}),
(\ref{eq3}), (\ref{eq5}), (\ref{eq7})-(\ref{eq9}).

\begin{theorem} \label{the51} The following inequalities among the
divergence measures hold:
\begin{equation}
\label{eq131} \frac{1}{4}\Delta (P\vert \vert Q) \leqslant
I(P\vert \vert Q) \leqslant h(P\vert \vert Q) \leqslant
\frac{1}{8}J(P\vert \vert Q) \leqslant T(P\vert \vert Q) \leqslant
\frac{1}{16}\Psi (P\vert \vert Q).
\end{equation}
\end{theorem}

The proof of the above theorem is based on the following
propositions, where we have proved each part separately.

\begin{proposition} The following inequality hold:
\begin{equation}
\label{eq132} \frac{1}{4}\Delta (P\vert \vert Q) \leqslant
I(P\vert \vert Q).
\end{equation}
\end{proposition}

\begin{proof} Let us consider
\begin{equation}
\label{eq133} g_{I\Delta } (x) = \frac{{f}''_I (x)}{{f}''_\Delta
(x)} = \frac{(x + 1)^2}{16x}, \quad x \in (0,\infty ),
\end{equation}

\noindent where ${f}''_I (x)$ and ${f}''_\Delta (x)$ are as given
by (\ref{eq37}) and (\ref{eq31}) respectively.

From (\ref{eq133}), we have
\begin{equation}
\label{eq134} {g}'_{I\Delta } (x) = \frac{(x - 1)(x + 1)}{16x^2}
\begin{cases}
 { \geqslant 0,} & {x \geqslant 1} \\
 { \leqslant 0,} & {x \leqslant 1} \\
\end{cases}.
\end{equation}

In view of (\ref{eq134}), we conclude that the function
$g_{I\Delta } (x)$ is decreasing in $x \in (0,1)$ and increasing
in $x \in (1,\infty )$, and hence
\begin{equation}
\label{eq135} m = \mathop {\sup }\limits_{x \in (0,\infty )}
g_{I\Delta } (x) = g_{I\Delta } (1) = \frac{1}{4}.
\end{equation}

Applying the inequalities (\ref{eq21}) for the measures $\Delta
(P\vert \vert Q)$ and $I(P\vert \vert Q)$ along with (\ref{eq135})
we get the required result.
\end{proof}

\begin{proposition} The following inequality hold:
\begin{equation}
\label{eq136} I(P\vert \vert Q) \leqslant h(P\vert \vert Q).
\end{equation}
\end{proposition}

\begin{proof} Let us consider
\begin{equation}
\label{eq137} g_{Ih} (x) = \frac{{f}''_I (x)}{{f}''_h (x)} =
\frac{2\sqrt x }{x + 1}, \quad x \in (0,\infty ),
\end{equation}

\noindent where ${f}''_I (x)$ and ${f}''_h (x)$ are as given by
(\ref{eq37}) and (\ref{eq29}) respectively.

From (\ref{eq137}), we have
\begin{equation}
\label{eq138} {g}'_{Ih} (x) = - \frac{x - 1}{\sqrt x (x + 1)^2}
\begin{cases}
 { \geqslant 0,} & {x \leqslant 1} \\
 { \leqslant 0,} & {x \geqslant 1} \\
\end{cases}.
\end{equation}

In view of (\ref{eq138}), we conclude that the function $g_{Ih}
(x)$ is increasing in $x \in (0,1)$ and decreasing in $x \in
(1,\infty )$, and hence
\begin{equation}
\label{eq139} M = \mathop {\sup }\limits_{x \in (0,\infty )}
g_{Ih} (x) = g_{Ih} (1) = 1.
\end{equation}

Applying the inequalities (\ref{eq21}) for the measures $I(P\vert
\vert Q)$ and $h(P\vert \vert Q)$ along with (\ref{eq139}) we get
the required result.
\end{proof}

\begin{proposition} The following inequality hold:
\begin{equation}
\label{eq140} h(P\vert \vert Q) \leqslant \frac{1}{8}J(P\vert
\vert Q).
\end{equation}
\end{proposition}

\begin{proof} Let us consider
\begin{equation}
\label{eq141} g_{Jh} (x) = \frac{{f}''_J (x)}{{f}''_h (x)} =
\frac{4(x + 1)}{\sqrt x }, \quad x \in (0,\infty ),
\end{equation}

\noindent where ${f}''_J (x)$ and ${f}''_h (x)$ are as given by
(\ref{eq35}) and (\ref{eq29}) respectively.

From (\ref{eq141}) we have
\begin{equation}
\label{eq142} {g}'_{Jh} (x) = \frac{2(x - 1)}{x\sqrt x }
\begin{cases}
 { \geqslant 0,} & {x \geqslant 1} \\
 { \leqslant 0,} & {x \leqslant 1} \\
\end{cases}.
\end{equation}

In view of (\ref{eq142}), we conclude that the function $g_{Jh}
(x)$ is decreasing in $x \in (0,1)$ and increasing in $x \in
(1,\infty )$, and hence
\begin{equation}
\label{eq143} m = \mathop {\inf }\limits_{x \in (0,\infty )}
g_{Jh} (x) = g_{Jh} (1) = 8.
\end{equation}

Applying the inequalities (\ref{eq21}) for the measures $h(P\vert
\vert Q)$ and $J(P\vert \vert Q)$ along with (\ref{eq143}) we get
the required result.
\end{proof}

\begin{proposition} The following inequality hold:
\begin{equation}
\label{eq144} \frac{1}{8}J(P\vert \vert Q) \leqslant T(P\vert
\vert Q).
\end{equation}
\end{proposition}

\begin{proof} Let us consider
\begin{equation}
\label{eq145} g_{JT} (x) = \frac{{f}''_J (x)}{{f}''_T (x)} =
\frac{4(x + 1)^2}{x^2 + 1}, \quad x \in (0,\infty ),
\end{equation}

\noindent where ${f}''_J (x)$ and ${f}''_T (x)$ are as given by
(\ref{eq35}) and (\ref{eq39}) respectively.

From (\ref{eq145}) we have
\begin{equation}
\label{eq146} {g}'_{JT} (x) = - \frac{8(x - 1)(x + 1)}{(x^2 +
1)^2} \begin{cases}
 { \geqslant 0,} & {x \leqslant 1} \\
 { \leqslant 0,} & {x \geqslant 1} \\
\end{cases}.
\end{equation}

In view of (\ref{eq146}) we conclude that the function $g_{JT}
(x)$ is increasing in $x \in (0,1)$ and decreasing in $x \in
(1,\infty )$, and hence
\begin{equation}
\label{eq147} M = \mathop {\sup }\limits_{x \in (0,\infty )}
g_{JT} (x) = g_{JT} (1) = 8.
\end{equation}

Applying the inequality (\ref{eq21}) for the measures $J(P\vert
\vert Q)$ and $T(P\vert \vert Q)$ along with (\ref{eq147}) we get
the required result.
\end{proof}

\begin{proposition} The following inequality hold:
\begin{equation}
\label{eq148} T(P\vert \vert Q) \leqslant \frac{1}{16}\Psi (P\vert
\vert Q).
\end{equation}
\end{proposition}

\begin{proof} Let us consider
\begin{equation}
\label{eq149} g_{T\Psi } (x) = \frac{{f}''_T (x)}{{f}''_\Psi (x)}
= \frac{x(x^2 + 1)}{8(x + 1)(x^3 + 1)}, \quad x \in (0,\infty ),
\end{equation}

\noindent where ${f}''_T (x)$ and ${f}''_\Psi (x)$ are as given by
(\ref{eq39}) and (\ref{eq33}) respectively.

From (\ref{eq149}) we have
\begin{equation}
\label{eq150} {g}'_{T\Psi } (x) = - \frac{(x - 1)(x^4 + 4x^2 +
1)}{8(x + 1)^3(x^2 - x + 1)^2} \begin{cases}
 { \geqslant 0,} & {x \leqslant 1} \\
 { \leqslant 0,} & {x \geqslant 1} \\
\end{cases}.
\end{equation}

In view of (\ref{eq150}) we conclude that the function $g_{T\Psi }
(x)$ is increasing in $x \in (0,1)$ and decreasing in $x \in
(1,\infty )$, and hence
\begin{equation}
\label{eq151} M = \mathop {\sup }\limits_{x \in (0,\infty )}
g_{T\Psi } (x) = g_{T\Psi } (1) = \frac{1}{16}.
\end{equation}

Applying the inequality (\ref{eq21}) for the measures $T(P\vert
\vert Q)$ and $\Psi (P\vert \vert Q)$ along with (\ref{eq151}) we
get the required result.
\end{proof}

The proof of the inequalities given in (\ref{eq131}) follows by
combining the results given in (\ref{eq132}), (\ref{eq136}),
(\ref{eq140}), 5.14) and (\ref{eq148}) respectively.\\

Dragomir et al. \cite{dsb} proved the following two inequalities
involving the measures (\ref{eq3}), (\ref{eq5}) and (\ref{eq7}):
\begin{equation}
\label{eq152} 0 \leqslant \frac{1}{2}J(P\vert \vert Q) - \Delta
(P\vert \vert Q) \leqslant \frac{1}{12}D^ * (P\vert \vert Q),
\end{equation}

\noindent and
\begin{equation}
\label{eq153} 0 \leqslant \frac{1}{2}\Psi (P\vert \vert Q) -
J(P\vert \vert Q) \leqslant \frac{1}{6}D^ * (P\vert \vert Q),
\end{equation}

\noindent where
\begin{equation}
\label{eq154} D^ * (P\vert \vert Q) = \sum\limits_{i = 1}^n
{\frac{(p_i - q_i )^4}{\sqrt {(p_i q_i )^3} }} .
\end{equation}

In the following section we shall improve the inequalities given in
(\ref{eq131}). An improvement over the inequalities (\ref{eq152})
and (\ref{eq153}) along with their unification is also presented.

\section{Difference of Divergence Measures}

Let us consider the following \textit{nonnegative} differences:
\begin{align}
\label{eq40} D_{\Psi T} (P\vert \vert Q) & = \frac{1}{16}\Psi
(P\vert \vert Q) - T(P\vert \vert Q), \\
\label{eq41} D_{\Psi J} (P\vert \vert Q) & = \frac{1}{16}\Psi
(P\vert \vert Q) - \frac{1}{8}J(P\vert \vert Q), \\
\label{eq42} D_{\Psi h} (P\vert \vert Q) & = \frac{1}{16}\Psi
(P\vert \vert Q) - h(P\vert \vert Q),\\
\label{eq43} D_{\Psi I} (P\vert \vert Q) & = \frac{1}{16}\Psi
(P\vert \vert Q) - I(P\vert \vert Q), \\
\label{eq44} D_{\Psi \Delta
} (P\vert \vert Q) & = \frac{1}{16}\Psi (P\vert \vert Q) -
\frac{1}{4}\Delta (P\vert
\vert Q), \\
\label{eq45} D_{TJ} (P\vert \vert Q) & = T(P\vert \vert Q) -
\frac{1}{8}J(P\vert \vert Q),\\
\label{eq46} D_{Th} (P\vert \vert Q) & = T(P\vert \vert Q) -
h(P\vert \vert Q),\\
\label{eq47} D_{TI} (P\vert \vert Q) & = T(P\vert \vert Q) -
I(P\vert \vert Q), \\
\label{eq48} D_{T\Delta } (P\vert \vert Q) & =
T(P\vert \vert Q) -
\frac{1}{4}\Delta (P\vert \vert Q),\\
\label{eq49} D_{Jh} (P\vert \vert Q) & = \frac{1}{8}J(P\vert \vert
Q) - h(P\vert \vert Q),\\
\label{eq50} D_{JI} (P\vert \vert Q) & = \frac{1}{8}J(P\vert \vert
Q) - I(P\vert \vert Q),\\
\label{eq51} D_{J\Delta } (P\vert \vert Q) & = \frac{1}{8}J(P\vert
\vert Q) - \frac{1}{4}\Delta (P\vert \vert Q),\\
\label{eq52} D_{hI} (P\vert \vert Q) & = h(P\vert \vert Q) -
I(P\vert \vert Q),\\
\label{eq53} D_{h\Delta } (P\vert \vert Q) & = h(P\vert \vert Q) -
\frac{1}{4}\Delta (P\vert \vert Q),\\
\intertext{and} \label{eq54} D_{I\Delta } (P\vert \vert Q) & =
I(P\vert \vert Q) - \frac{1}{4}\Delta (P\vert \vert Q).
\end{align}

In the examples below we shall show the convexity of the above
measures (\ref{eq40})-(\ref{eq54}). In view of Theorem 2.1 and
Examples 2.1-2.6, it is sufficient to show the nonnegativity of
the second order derivative of generating function in each case.

\begin{example} \label{exa31} We can write
\[
D_{\Psi T} (P\vert \vert Q) = \frac{1}{16}\Psi (P\vert \vert Q) - T(P\vert
\vert Q)
 = \sum\limits_{i = 1}^n {q_i f_{\Psi T} \left( {\frac{p_i }{q_i }} \right)}
,
\]

\noindent where
\[
f_{\Psi T} (x) = \frac{1}{16}f_\Psi \left( x \right) - f_T \left(
x \right), \, x > 0.
\]

Moreover, we have
\begin{align}
\label{eq55} {f}''_{\Psi T} (x) & = \frac{1}{16}{f}''_\Psi \left(
x \right) - {f}''_T \left( x \right)\\
  & = \frac{x^3 + 1}{8x^3} - \frac{x^2 + 1}{4x^2(x + 1)}
= \frac{(x - 1)^2(x^2 + x + 1)}{8x^3(x + 1)}
 \geqslant 0,
\, \forall x > 0,\notag
\end{align}

\noindent where ${f}''_\Psi \left( x \right)$ and ${f}''_T \left(
x \right)$ are as given by (\ref{eq33}) and (\ref{eq39})
respectively.
\end{example}

\begin{example} \label{exa32} We can write
\[
D_{\Psi J} (P\vert \vert Q) = \frac{1}{16}\Psi (P\vert \vert Q) -
\frac{1}{8}J(P\vert \vert Q)
 = \sum\limits_{i = 1}^n {q_i f_{\Psi J} \left( {\frac{p_i }{q_i }} \right)}
,
\]

\noindent where
\[
f_{\Psi J} (x) = \frac{1}{16}f_\Psi \left( x \right) -
\frac{1}{8}f_J \left( x \right), \, x > 0.
\]

Moreover, we have
\begin{align}
\label{eq56} {f}''_{\Psi J} (x) & = \frac{1}{16}{f}''_\Psi \left(
x \right) - \frac{1}{8}{f}''_J \left( x \right) \\
 & = \frac{1}{8}\left( {\frac{x^3 + 1}{x^3} - \frac{x + 1}{x^2}}
 \right)
= \frac{(x - 1)^2(x + 1)}{8x^3}
 \geqslant 0, \, \forall x > 0,\notag
\end{align}

\noindent where ${f}''_\Psi \left( x \right)$ and ${f}''_J \left(
x \right)$ are as given by (\ref{eq33}) and (\ref{eq35})
respectively.
\end{example}

\begin{example} \label{exa33} We can write
\[
D_{\Psi h} (P\vert \vert Q) = \frac{1}{16}\Psi (P\vert \vert Q) - h(P\vert
\vert Q)
 = \sum\limits_{i = 1}^n {q_i f_{\Psi h} \left( {\frac{p_i }{q_i }} \right)}
,
\]

\noindent where
\[
f_{\Psi h} (x) = \frac{1}{16}f_\Psi \left( x \right) - f_h \left(
x \right), \, x > 0.
\]

Moreover, we have
\begin{align}
\label{eq57} {f}''_{\Psi h} (x) & = \frac{1}{16}{f}''_\Psi \left(
x \right) - {f}''_h \left( x \right)\\
 & = \frac{1}{4}\left( {\frac{x^3 + 1}{2x^3} - \frac{1}{x\sqrt x }}
 \right)
 = \frac{\left( {x\sqrt x - 1} \right)^2}{8x^3}
 \geqslant 0,
\, \forall x > 0,\notag
\end{align}

\noindent where ${f}''_\Psi \left( x \right)$ and ${f}''_h \left(
x \right)$ are as given by (\ref{eq33}) and (\ref{eq29})
respectively.
\end{example}

\begin{example} \label{exa34} We can write
\[
D_{\Psi I} (P\vert \vert Q) = \frac{1}{16}\Psi (P\vert \vert Q) -
I(P\vert \vert Q)
 = \sum\limits_{i = 1}^n {q_i f_{\Psi I} \left( {\frac{p_i }{q_i }} \right)}
,
\]

\noindent where
\[
f_{\Psi I} (x) = \frac{1}{16}f_\Psi \left( x \right) - f_I \left(
x \right), \, x > 0.
\]

Moreover, we have
\begin{align}
\label{eq58} {f}''_{\Psi I} (x) & = \frac{1}{16}{f}''_\Psi \left(
x \right) - {f}''_I \left( x \right)\\
 & = \frac{1}{2x}\left( {\frac{x^3 + 1}{4x^2} - \frac{1}{x + 1}}
 \right)
 = \frac{(x - 1)^2(x^2 + 3x + 1)}{8x^3(x + 1)}
 \geqslant 0,
\, \forall x > 0,\notag
\end{align}

\noindent where ${f}''_\Psi \left( x \right)$ and ${f}''_I \left(
x \right)$ are as given by (\ref{eq33}) and (\ref{eq37})
respectively.
\end{example}

\begin{example} \label{exa35} We can write
\[
D_{\Psi \Delta } (P\vert \vert Q) = \frac{1}{16}\Psi (P\vert \vert Q) -
\frac{1}{4}\Delta (P\vert \vert Q)
 = \sum\limits_{i = 1}^n {q_i f_{\Psi \Delta } \left( {\frac{p_i }{q_i }}
\right)} ,
\]

\noindent where
\[
f_{\Psi \Delta } (x) = \frac{1}{4}\left( {\frac{1}{4}f_\Psi \left(
x \right) - f_\Delta \left( x \right)} \right), \, x > 0.
\]

Moreover, we have
\begin{align}
\label{eq59} {f}''_{\Psi \Delta } (x) & = \frac{1}{4}\left(
{\frac{1}{4}{f}''_\Psi \left( x \right) - {f}''_\Delta \left( x
\right)} \right)
= \frac{x^3 + 1}{8x^3} - \frac{2}{(x + 1)^3}\\
 & = \frac{(x - 1)^2(x^4 + 5x^3 + 12x^2 + 5x + 1)}{8x^3(x + 1)^3}
 \geqslant 0,
\, \forall x > 0,\notag
\end{align}

\noindent where ${f}''_\Psi \left( x \right)$ and ${f}''_\Delta
\left( x \right)$ are as given by (\ref{eq33}) and (\ref{eq31})
respectively.
\end{example}

\begin{example} \label{exa36} We can write
\[
D_{TJ} (P\vert \vert Q) = T(P\vert \vert Q) - \frac{1}{8}J(P\vert \vert Q)
 = \sum\limits_{i = 1}^n {q_i f_{TJ} \left( {\frac{p_i }{q_i }} \right)} ,
\]

\noindent where
\[
f_{TJ} (x) = f_T \left( x \right) - \frac{1}{8}f_J \left( x
\right), \, x > 0.
\]

Moreover, we have
\begin{align}
\label{eq60} {f}''_{TJ} (x) & = {f}''_T \left( x \right) -
\frac{1}{8}{f}''_J \left( x \right) \\
& = \frac{x^2+1}{4x^2(x + 1)} - \frac{x + 1}{8x^2} = \frac{(x -
1)^2}{8x^2(x + 1)}
 \geqslant 0,
\, \forall x > 0,\notag
\end{align}

\noindent where ${f}''_T \left( x \right)$ and ${f}''_J \left( x
\right)$ are as given by (\ref{eq39}) and (\ref{eq35})
respectively.
\end{example}

\begin{example} \label{exa37} We can write
\[
D_{Th} (P\vert \vert Q) = T(P\vert \vert Q) - h(P\vert \vert Q)
 = \sum\limits_{i = 1}^n {q_i f_{Th} \left( {\frac{p_i }{q_i }} \right)} ,
\]

\noindent where
\[
f_{Th} (x) = f_T \left( x \right) - f_h \left( x \right), \, x >
0.
\]

Moreover, we have
\begin{align}
\label{eq61} {f}''_{Th} (x) & = {f}''_T \left( x \right) - {f}''_h
\left( x \right)
 = \frac{1}{4}\left( {\frac{x^2+1}{x^2(x+1)} - \frac{1}{x\sqrt x
 }}\right)\\
 & = \frac{\left( {\sqrt x - 1} \right)^2\left( {x + \sqrt x + 1}
\right)}{4x^2(x + 1)}
 \geqslant 0,
\, \forall x > 0,\notag
\end{align}

\noindent where ${f}''_T \left( x \right)$ and ${f}''_h \left( x
\right)$ are as given by (\ref{eq39}) and (\ref{eq29})
respectively.
\end{example}

\begin{example} \label{38} We can write
\[
D_{TI} (P\vert \vert Q) = T(P\vert \vert Q) - I(P\vert \vert Q)
 = \sum\limits_{i = 1}^n {q_i f_{TI} \left( {\frac{p_i }{q_i }} \right)} ,
\]

\noindent where
\[
f_{TI} (x) = f_T \left( x \right) - f_I \left( x \right), \, x >
0.
\]

Moreover, we have
\begin{align}
\label{eq62} {f}''_{TI} (x) & = {f}''_T \left( x \right) - {f}''_I
\left( x \right)\\
 & = \frac{x^2+1}{4x^2(x + 1)} - \frac{1}{2x(x + 1)}
 = \frac{(x - 1)^2}{4x^2(x + 1)}
 \geqslant 0,
\, \forall x > 0,\notag
\end{align}

\noindent where ${f}''_T \left( x \right)$ and ${f}''_I \left( x
\right)$ are as given by (\ref{eq39}) and (\ref{eq37})
respectively.
\end{example}

\begin{example} \label{exa39} We can write
\[
D_{T\Delta } (P\vert \vert Q) = T(P\vert \vert Q) - \frac{1}{4}\Delta
(P\vert \vert Q)
 = \sum\limits_{i = 1}^n {q_i f_{T\Delta } \left( {\frac{p_i }{q_i }}
\right)} ,
\]

\noindent where
\[
f_{T\Delta } (x) = f_T \left( x \right) - \frac{1}{4}f_\Delta
\left( x \right), \, x > 0.
\]

Moreover, we have
\begin{align}
\label{eq63} {f}''_{T\Delta } (x) & = {f}''_T \left( x \right) -
\frac{1}{4}{f}''_\Delta \left( x \right)
 = \frac{x^2 + 1}{4x^2(x + 1)} - \frac{8}{(x + 1)^3}\\
 & = \frac{(x - 1)^2(x^2 + 4x + 1)}{4x^2(x + 1)^3}
 \geqslant 0,
\, \forall x > 0,\notag
\end{align}

\noindent where ${f}''_T \left( x \right)$ and ${f}''_\Delta
\left( x \right)$ are as given by (\ref{eq39}) and (\ref{eq31})
respectively.
\end{example}

\begin{example} \label{310} We can write
\[
D_{Jh} (P\vert \vert Q) = \frac{1}{8}J(P\vert \vert Q) - h(P\vert \vert Q)
 = \sum\limits_{i = 1}^n {q_i f_{Jh} \left( {\frac{p_i }{q_i }} \right)} ,
\]

\noindent where
\[
f_{Jh} (x) = \frac{1}{8}f_J \left( x \right) - f_h \left( x
\right), \, x > 0.
\]

Moreover, we have
\begin{align}
\label{eq64} {f}''_{Jh} (x) & = \frac{1}{8}{f}''_J \left( x
\right) - {f}''_h \left( x \right)\\
 & = \frac{x + 1}{8x^2} - \frac{1}{4x\sqrt x }
 = \frac{\left( {\sqrt x - 1} \right)^2}{8x^2}
 \geqslant 0,
\, \forall x > 0,\notag
\end{align}

\noindent where ${f}''_J \left( x \right)$ and ${f}''_h \left( x
\right)$ are as given by (\ref{eq35}) and (\ref{eq29})
respectively.
\end{example}

\begin{example} \label{exa311} We can write
\[
D_{JI} (P\vert \vert Q) = \frac{1}{8}J(P\vert \vert Q) - I(P\vert \vert Q)
 = \sum\limits_{i = 1}^n {q_i f_{JI} \left( {\frac{p_i }{q_i }} \right)} ,
\]

\noindent where
\[
f_{JI} (x) = \frac{1}{8}f_J \left( x \right) - f_I \left( x
\right), \, x > 0.
\]

Moreover, we have
\begin{align}
\label{eq65} {f}''_{JI} (x) & = \frac{1}{8}{f}''_J \left( x
\right) - {f}''_I \left( x \right) \\
 & = \frac{x + 1}{8x^2} - \frac{1}{2x(x + 1)}
 = \frac{(x - 1)^2}{8x^2(x + 1)}
 \geqslant 0,
\, \forall x > 0,\notag
\end{align}

\noindent where ${f}''_J \left( x \right)$ and ${f}''_I \left( x
\right)$ are as given by (\ref{eq35}) and (\ref{eq37})
respectively.
\end{example}

\begin{example} \label{312} We can write
\[
D_{J\Delta } (P\vert \vert Q) = \frac{1}{8}J(P\vert \vert Q) -
\frac{1}{4}\Delta (P\vert \vert Q)
 = \sum\limits_{i = 1}^n {q_i f_{J\Delta } \left( {\frac{p_i }{q_i }}
\right)} ,
\]

\noindent where
\[
f_{J\Delta } (x) = \frac{1}{8}f_J \left( x \right) -
\frac{1}{4}f_\Delta \left( x \right), \, x > 0.
\]

Moreover, we have
\begin{align}
\label{eq66} {f}''_{J\Delta } (x) & = \frac{1}{8}{f}''_J \left( x
\right) - \frac{1}{4}{f}''_\Delta \left( x \right)
 = \frac{x + 1}{8x^2} - \frac{2}{(x + 1)^3}\\
 & = \frac{(x - 1)^2(x^2 + 6x + 1)}{8x^2(x + 1)^3}
 \geqslant 0, \, \forall x > 0,\notag
\end{align}

\noindent where ${f}''_J \left( x \right)$ and ${f}''_\Delta
\left( x \right)$ are as given by (\ref{eq35}) and (\ref{eq31})
respectively.
\end{example}

\begin{example} \label{313} We can write
\[
D_{hI} (P\vert \vert Q) = h(P\vert \vert Q) - I(P\vert \vert Q)
 = \sum\limits_{i = 1}^n {q_i f_{hI} \left( {\frac{p_i }{q_i }} \right)} ,
\]

\noindent where
\[
f_{hI} (x) = f_h \left( x \right) - f_I \left( x \right), \, x >
0.
\]

Moreover, we have
\begin{align}
\label{eq67} {f}''_{hI} (x) & = {f}''_h \left( x \right) - {f}''_I
\left( x \right)\\
&  = \frac{1}{4x\sqrt x } - \frac{1}{2x(x + 1)}
 = \frac{(\sqrt x - 1)^2}{4x^{3 / 2}(x + 1)}
 \geqslant 0, \, \forall x > 0,\notag
\end{align}

\noindent where ${f}''_h \left( x \right)$ and ${f}''_I \left( x
\right)$ are as given by (\ref{eq29}) and (\ref{eq37})
respectively.
\end{example}

\begin{example} \label{exa314} We can write
\[
D_{h\Delta } (P\vert \vert Q) = h(P\vert \vert Q) - \frac{1}{4}\Delta
(P\vert \vert Q)
 = \sum\limits_{i = 1}^n {q_i f_{h\Delta } \left( {\frac{p_i }{q_i }}
\right)} ,
\]

\noindent where
\[
f_{h\Delta } (x) = f_h \left( x \right) - \frac{1}{4}f_\Delta
\left( x \right), \, x > 0.
\]

Moreover, we have
\begin{align}
\label{eq68} {f}''_{h\Delta } (x) & = {f}''_h \left( x \right) -
\frac{1}{4}{f}''_\Delta \left( x \right)
 = \frac{1}{4x\sqrt x } - \frac{2}{(x + 1)^3}\\
 & = \frac{\left( {\sqrt x - 1} \right)^2\left[ {\left( {\sqrt x + 1}
\right)^2\left( {x + 1} \right) + 4x} \right]}{4x^{3 / 2}(x + 1)^3}
 \geqslant 0, \, \forall x > 0,\notag
\end{align}

\noindent where ${f}''_h \left( x \right)$ and ${f}''_\Delta
\left( x \right)$ are as given by (\ref{eq29}) and (\ref{eq31})
respectively.
\end{example}

\begin{example} \label{exa315} We can write
\[
D_{I\Delta } (P\vert \vert Q) = I(P\vert \vert Q) - \frac{1}{4}\Delta
(P\vert \vert Q)
 = \sum\limits_{i = 1}^n {q_i f_{I\Delta } \left( {\frac{p_i }{q_i }}
\right)} ,
\]

\noindent where
\[
f_{I\Delta } (x) = f_I \left( x \right) - \frac{1}{4}f_\Delta
\left( x \right), \, x > 0.
\]

Moreover, we have
\begin{align}
\label{eq69} {f}''_{I\Delta } (x) & = {f}''_I \left( x \right) -
\frac{1}{4}{f}''_\Delta \left( x \right)\\
& = \frac{1}{2x(x + 1)} - \frac{2}{(x + 1)^3}
 = \frac{(x - 1)^2}{2x(x + 1)^3}
 \geqslant 0, \, \forall x > 0,\notag
\end{align}

\noindent where ${f}''_I \left( x \right)$ and ${f}''_\Delta
\left( x \right)$ are as given by (\ref{eq37}) and (\ref{eq31})
respectively.
\end{example}

Thus in view of Theorem \ref{the21} and Examples
\ref{exa31}-\ref{exa315}, we can say that the \textit{divergence
measures} given in (\ref{eq40})-(\ref{eq54}) are all
\textit{nonnegative} and \textit{convex} in the pair of
probability distributions $(P,Q) \in \Gamma _n \times \Gamma _n $.

\section{Refinement Inequalities}

In view of (\ref{eq131}), the following inequalities are obvious:
\begin{align}
\label{eq70} & D_{\Psi T} (P\vert \vert Q) \leqslant D_{\Psi J}
(P\vert \vert Q) \leqslant D_{\Psi h} (P\vert \vert Q) \leqslant
D_{\Psi I} (P\vert \vert Q) \leqslant D_{\Psi \Delta } (P\vert
\vert Q), \\
\label{eq71} & D_{TJ} (P\vert \vert Q) \leqslant D_{Th} (P\vert
\vert Q) \leqslant D_{TI} (P\vert \vert Q) \leqslant D_{T\Delta }
(P\vert \vert Q), \\
\label{eq72} & D_{Jh} (P\vert \vert Q) \leqslant D_{JI} (P\vert
\vert Q) \leqslant D_{J\Delta } (P\vert \vert Q)\\
\intertext{and} \label{eq73} & D_{hI} (P\vert \vert Q) \leqslant
D_{h\Delta } (P\vert \vert Q).
\end{align}

\bigskip
In view of the relation (\ref{eq10}), we have the following
equality:
\begin{equation}
\label{eq77} D_{JI} (P\vert \vert Q) = \frac{1}{2}D_{TI} (P\vert
\vert Q) = D_{TJ} (P\vert \vert Q).
\end{equation}

\bigskip
In this section our aim is to establish refinement inequalities
improving the one given in (\ref{eq131}). This refinement is given
in the following theorem.

\begin{theorem} \label{the41} The following inequalities hold:
\begin{align}
\label{eq74} & D_{I\Delta } (P\vert \vert Q) \leqslant
\frac{2}{3}D_{h\Delta } (P\vert \vert Q) \leqslant 2D_{hI} (P\vert
\vert Q) \leqslant D_{TJ} (P\vert \vert Q), \\
\label{eq75} & D_{I\Delta } (P\vert \vert Q) \leqslant
\frac{2}{3}D_{h\Delta } (P\vert \vert Q) \leqslant
\frac{1}{2}D_{J\Delta } (P\vert \vert Q) \leqslant
\frac{1}{3}D_{T\Delta } (P\vert \vert Q) \leqslant D_{TJ} (P\vert
\vert Q),\\
\intertext{and}
\label{eq76} & D_{TJ} (P\vert \vert Q) \leqslant
\frac{2}{3}D_{Th} (P\vert \vert Q) \leqslant 2D_{Jh} (P\vert \vert
Q) \leqslant
\frac{1}{6}D_{\Psi \Delta } (P\vert \vert Q)\\
 & \qquad \leqslant \frac{1}{5}D_{\Psi I} (P\vert \vert Q) \leqslant
\frac{2}{9}D_{\Psi h} (P\vert \vert Q) \leqslant \frac{1}{4}D_{\Psi
J} (P\vert \vert Q) \leqslant \frac{1}{3}D_{\Psi T} (P\vert \vert
Q), \notag
\end{align}
\end{theorem}

The proofs of the inequalities (\ref{eq74})-(\ref{eq76}) are based
on the following propositions.

\begin{proposition} \label{pro41} We have
\begin{equation}
\label{eq78}
D_{I\Delta } (P\vert \vert Q) \leqslant \frac{2}{3}D_{h\Delta } (P\vert
\vert Q).
\end{equation}
\end{proposition}

\begin{proof} Let us consider
\begin{align}
g_{I\Delta \_h\Delta } (x) = \frac{{f}''_{I\Delta }
(x)}{{f}''_{h\Delta } (x)} & = \frac{2\sqrt x \left( {x - 1}
\right)^2}{\left( {x + 1} \right)^3 - 8\left( {\sqrt x }
\right)^3}, \, x \ne 1\notag\\
 & = \frac{2\sqrt x \left( {\sqrt x + 1} \right)^2}{\left( {\sqrt x + 1}
\right)^2(x + 1) + 4x}\notag
\end{align}

\noindent for all $x \in (0,\infty )$, where ${f}''_{I\Delta }
(x)$ and ${f}''_{h\Delta } (x)$ are as given by (\ref{eq69}) and
(\ref{eq68}) respectively.

Calculating the first order derivative of the function $g_{I\Delta
\_h\Delta } (x)$ with respect to $x$, one gets
\begin{align}
\label{eq79} {g}'_{I\Delta \_h\Delta } (x) & = - \frac{\left(
{\sqrt x + 1} \right)\left( {x^{5 / 2} - 2x^{3 / 2} + 3x^2 + 2x -
3\sqrt x - 1} \right)}{\sqrt x \left[ {x^2 + 6x + 2\sqrt x \left(
{x + 1} \right) + 1} \right]^2}\\
& = - \frac{(x - 1)(x + 1)\left( {x + 4\sqrt x + 1} \right)}{\sqrt
x \left[ {x^2 + 6x + 2\sqrt x \left( {x + 1} \right) + 1}
\right]^2}
\begin{cases}
 { > 0,} & {x < 1} \\
 { < 0,} & {x > 1} \\
\end{cases}.\notag
\end{align}

In view of (\ref{eq79}) we conclude that the function $g_{I\Delta
\_h\Delta } (x)$ is increasing in $x \in (0,1)$ and decreasing in
$x \in (1,\infty )$, and hence
\begin{equation}
\label{eq80}
M = \mathop {\sup }\limits_{x \in (0,\infty )} g_{I\Delta \_h\Delta } (x) =
g_{I\Delta \_h\Delta } (1) = \frac{2}{3}.
\end{equation}

By the application of (\ref{eq21}) with (\ref{eq80}) we get
(\ref{eq78}).
\end{proof}

\begin{proposition} \label{pro42} We have
\begin{equation}
\label{eq81}
D_{h\Delta } (P\vert \vert Q) \leqslant 3D_{hI} (P\vert \vert Q).
\end{equation}
\end{proposition}

\begin{proof} Let us consider
\[
g_{h\Delta \_hI} (x) = \frac{{f}''_{h\Delta } (x)}{{f}''_{hI} (x)}
= \frac{(x + 1)\left( {\sqrt x + 1} \right)^2 + 4x}{(x + 1)^2}, \,
x \in (0,\infty ),
\]

\noindent where ${f}''_{h\Delta } (x)$ and ${f}''_{hI} (x)$ are as
given by (\ref{eq68}) and (\ref{eq67}) respectively.

Calculating the first order derivative of the function $g_{h\Delta
\_hI} (x)$ with respect to $x$, one gets
\begin{align}
\label{eq82} {g}'_{h\Delta \_hI} (x) & = - \frac{4x^{3 / 2} + x^2
- 4\sqrt x - 1}{\sqrt x \left( {x + 1} \right)^3}\\
& = - \frac{(x - 1)\left( {x + 4\sqrt x + 1} \right)}{\sqrt x
\left( {x + 1} \right)^3} \,
\begin{cases}
 { > 0,} & {x < 1} \\
 { < 0,} & {x > 1} \\
\end{cases}.\notag
\end{align}

In view of (\ref{eq82}) we conclude that the function $g_{h\Delta
\_hI} (x)$ is increasing in $x \in (0,1)$ and decreasing in $x \in
(1,\infty )$, and hence
\begin{equation}
\label{eq83}
M = \mathop {\sup }\limits_{x \in (0,\infty )} g_{h\Delta \_hI} (x) =
g_{h\Delta \_hI} (1) = 3.
\end{equation}

By the application of (\ref{eq21}) with (\ref{eq83}) we get
(\ref{eq81}).
\end{proof}

\begin{remark} \label{not41} In view of Propositions \ref{pro41} and \ref{pro42},
and the inequality (\ref{eq131}) we conclude that
\begin{equation}
\label{eq84}
I(P\vert \vert Q) \leqslant \frac{2}{3}h(P\vert \vert Q) +
\frac{1}{12}\Delta (P\vert \vert Q) \leqslant h(P\vert \vert Q).
\end{equation}
\end{remark}

\begin{proposition} \label{pro43} We have
\begin{equation}
\label{eq85}
D_{hI} (P\vert \vert Q) \leqslant \frac{1}{2}D_{TJ} (P\vert \vert Q).
\end{equation}
\end{proposition}

\begin{proof} Let us consider
\[
g_{hI\_TJ} (x) = \frac{{f}''_{hI} (x)}{{f}''_{TJ} (x)} =
\frac{2\sqrt x }{\left( {\sqrt x + 1} \right)^2}, \, x \in
(0,\infty ),
\]

\noindent where ${f}''_{hI} (x)$ and ${f}''_{TJ} (x)$ are as given
by (\ref{eq67}) and (\ref{eq60}) respectively.

Calculating the first order derivative of the function $g_{hI\_TJ}
(x)$ with respect to $x$, one gets
\begin{equation}
\label{eq86} {g}'_{hI\_TJ} (x) = - \frac{\sqrt x - 1}{\sqrt x
\left( {\sqrt x + 1} \right)^3} \,
\begin{cases}
 { > 0,} & {x < 1} \\
 { < 0,} & {x > 1} \\
\end{cases}.
\end{equation}

In view of (\ref{eq86}), we conclude that the function $g_{hI\_TJ}
(x)$ is increasing in $x \in (0,1)$ and decreasing in $x \in
(1,\infty )$, and hence
\begin{equation}
\label{eq87}
M = \mathop {\sup }\limits_{x \in (0,\infty )} g_{hI\_TJ} (x) = g_{hI\_TJ}
(1) = \frac{1}{2}.
\end{equation}

By the application of (\ref{eq21}) with (\ref{eq87}) we get
(\ref{eq85}).
\end{proof}

\begin{remark} \label{not42} In view of Propositions \ref{pro43} and the inequality
(\ref{eq131}) we conclude the following inequality
\begin{equation}
\label{eq88}
h(P\vert \vert Q) \leqslant \frac{1}{16}J(P\vert \vert Q) +
\frac{1}{2}I(P\vert \vert Q) \leqslant \frac{1}{8}J(P\vert \vert Q).
\end{equation}
\end{remark}

Combining the inequalities (\ref{eq78}), (\ref{eq81}) and
(\ref{eq85}) we get (\ref{eq74}).

\begin{proposition} \label{pro44} We have
\begin{equation}
\label{eq89}
D_{h\Delta } (P\vert \vert Q) \leqslant \frac{3}{4}D_{J\Delta } (P\vert
\vert Q).
\end{equation}
\end{proposition}

\begin{proof} Let us consider
\begin{align}
g_{h\Delta \_J\Delta } (x) = \frac{{f}''_{h\Delta }
(x)}{{f}''_{J\Delta } (x)} & = \frac{2\sqrt x \left[ {(x + 1)^3 -
8x^{3 / 2}} \right]}{(x - 1)^2(x^2 + 6x + 1)}, \, x \ne 1\notag\\
 & = \frac{2\sqrt x \left[ {\left( {\sqrt x + 1} \right)^2(x + 1) + 4x}
\right]}{\left( {\sqrt x + 1} \right)^2\left( {x^2 + 6x + 1}
\right)},\notag
\end{align}

\noindent for all $x \in (0,\infty )$, where ${f}''_{h\Delta }
(x)$ and ${f}''_{J\Delta } (x)$ are as given by (\ref{eq68}) and
(\ref{eq66}) respectively.

Calculating the first order derivative of the function $g_{h\Delta
\_J\Delta } (x)$ with respect to $x$, one gets
\begin{align}
\label{eq90} {g}'_{h\Delta \_J\Delta } (x) & = - \frac{1}{\sqrt x
\left( {\sqrt x + 1} \right)^3\left( {x^2 + 6x + 1}
\right)^2}\left[ {3x^4 - 4x^3 - 18x^2 - 12x - 1} \right.\\
& \qquad \qquad \qquad \qquad \left. { + \sqrt x \left( {x^4 +
12x^3 + 18x^2 +
4x - 3} \right)} \right]\notag\\
& = - \frac{\left( {\sqrt x - 1} \right)(x+1)^2 \left( {x^2 +
4x\sqrt x + 14x + 4\sqrt x + 1} \right)}{\sqrt x \left( {\sqrt x +
1} \right)^3\left( {x^2 + 6x + 1} \right)^2} \,
\begin{cases}
 { > 0,} & {x < 1} \\
 { < 0,} & {x > 1} \\
\end{cases}.\notag
\end{align}

In view of (\ref{eq90}) we conclude that the function $g_{h\Delta
\_J\Delta } (x)$ is increasing in $x \in (0,1)$ and decreasing in
$x \in (1,\infty )$, and hence
\begin{equation}
\label{eq91}
M = \mathop {\sup }\limits_{x \in (0,\infty )} g_{h\Delta \_J\Delta } (x) =
g_{h\Delta \_J\Delta } (1) = \frac{3}{4}.
\end{equation}

By the application of (\ref{eq21}) with (\ref{eq91}) we get
(\ref{eq89}).
\end{proof}

\begin{remark} \label{not43} In view of Proposition \ref{pro44} and
the inequality (\ref{eq131}) we conclude the following inequality
\begin{equation}
\label{eq92}
h(P\vert \vert Q) \leqslant \frac{3}{32}J(P\vert \vert Q) +
\frac{1}{16}\Delta (P\vert \vert Q) \leqslant \frac{1}{8}J(P\vert \vert Q).
\end{equation}
\end{remark}

\begin{proposition} \label{pro45} We have
\begin{equation}
\label{eq93}
D_{J\Delta } (P\vert \vert Q) \leqslant \frac{2}{3}D_{T\Delta } (P\vert
\vert Q).
\end{equation}
\end{proposition}

\begin{proof} Let us consider
\[
g_{J\Delta \_T\Delta } (x) = \frac{{f}''_{J\Delta }
(x)}{{f}''_{T\Delta } (x)} = \frac{x^2 + 6x + 1}{2(x^2 + 4x + 1)},
\, x \in (0,\infty ),
\]

\noindent where ${f}''_{J\Delta } (x)$ and ${f}''_{T\Delta } (x)$
are as given by (\ref{eq66}) and (\ref{eq63}) respectively.

Calculating the first order derivative of the function $g_{J\Delta
\_T\Delta } (x)$ with respect to $x$, one gets
\begin{equation}
\label{eq94} {g}'_{J\Delta \_T\Delta } (x) = - \frac{(x - 1)(x +
1)}{(x^2 + 4x + 1)^2} \,
\begin{cases}
 { > 0,} & {x < 1} \\
 { < 0,} & {x > 1} \\
\end{cases}.
\end{equation}

In view of (\ref{eq94}) we conclude that the function $g_{J\Delta
\_T\Delta } (x)$ is increasing in $x \in (0,1)$ and decreasing in
$x \in (1,\infty )$, and hence
\begin{equation}
\label{eq95}
M = \mathop {\sup }\limits_{x \in (0,\infty )} g_{J\Delta \_T\Delta } (x) =
g_{J\Delta \_T\Delta } (1) = \frac{2}{3}.
\end{equation}

By the application of (\ref{eq21}) with (\ref{eq95}) we get
(\ref{eq93}).
\end{proof}

\begin{proposition} \label{pro46} We have
\begin{equation}
\label{eq96}
D_{T\Delta } (P\vert \vert Q) \leqslant 3D_{TJ} (P\vert \vert Q).
\end{equation}
\end{proposition}

\begin{proof} Let us consider
\[
g_{T\Delta \_TJ} (x) = \frac{{f}''_{T\Delta } (x)}{{f}''_{TJ} (x)}
= \frac{2(x^2 + 4x + 1)}{(x + 1)^2}, \, x \in (0,\infty ),
\]

\noindent where ${f}''_{T\Delta } (x)$ and ${f}''_{TJ} (x)$ are as
given by (\ref{eq63}) and (\ref{eq60}) respectively.

Calculating the first order derivative of the function $g_{T\Delta
\_TJ} (x)$ with respect to $x$, one gets
\begin{equation}
\label{eq97} {g}'_{T\Delta \_TJ} (x) = - \frac{4(x - 1)}{(x +
1)^3} \,
\begin{cases}
 { > 0,} & {x < 1} \\
 { < 0,} & {x > 1} \\
\end{cases}.
\end{equation}

In view of (\ref{eq97}) we conclude that the function $g_{T\Delta
\_TJ} (x)$ is increasing in $x \in (0,1)$ and decreasing in $x \in
(1,\infty )$, and hence
\begin{equation}
\label{eq98}
M = \mathop {\sup }\limits_{x \in (0,\infty )} g_{T\Delta \_TJ} (x) =
g_{T\Delta \_TJ} (1) = 3.
\end{equation}

By the application of (\ref{eq21}) with (\ref{eq98}) we get
(\ref{eq96}).
\end{proof}

\begin{remark} \label{not44} In view of Propositions \ref{pro45} and \ref{pro46},
and the inequality (\ref{eq131}) we conclude the following
inequality
\begin{equation}
\label{eq99}
\frac{1}{8}J(P\vert \vert Q) \leqslant \frac{2}{3}T(P\vert \vert Q) +
\frac{1}{12}\Delta (P\vert \vert Q) \leqslant T(P\vert \vert Q).
\end{equation}
\end{remark}

Combining the inequalities (\ref{eq78}), (\ref{eq89}), (\ref{eq93})
and (\ref{eq96}), we get (\ref{eq75}).

\begin{proposition} \label{pro47} We have
\begin{equation}
\label{eq100}
D_{TJ} (P\vert \vert Q) \leqslant \frac{2}{3}D_{Th} (P\vert \vert Q).
\end{equation}
\end{proposition}

\begin{proof} Let us consider
\begin{align}
g_{TJ\_Th} (x) = \frac{{f}''_{TJ} (x)}{{f}''_{Th} (x)} & =
\frac{(x - 1)^2}{2\left[ {x^2 + 1 - 2\sqrt x \left( {x + 1}
\right)} \right]}, \, x \ne 1 \,\notag\\
 & = \frac{\left( {\sqrt x + 1} \right)^2}{2\left( {x + \sqrt x + 1}
\right)},\notag
\end{align}

\noindent for all $x \in (0,\infty )$, where ${f}''_{TJ} (x)$ and
${f}''_{Th} (x)$ are as given by (\ref{eq60}) and (\ref{eq61})
respectively.

Calculating the first order derivative of the function $g_{TJ\_Th}
(x)$ with respect to $x$, one gets
\begin{equation}
\label{eq101} {g}'_{TJ\_Th} (x) = - \frac{\left( {\sqrt x - 1}
\right)\left( {\sqrt x + 1} \right)}{4\sqrt x \left( {x + \sqrt x
+ 1} \right)} \,
\begin{cases}
 { > 0,} & {x < 1} \\
 { < 0,} & {x > 1} \\
\end{cases}.
\end{equation}

In view of (\ref{eq101}) we conclude that the function $g_{TJ\_Th}
(x)$ is increasing in $x \in (0,1)$ and decreasing in $x \in
(1,\infty )$, and hence
\begin{equation}
\label{eq102}
M = \mathop {\sup }\limits_{x \in (0,\infty )} g_{TJ\_Th} (x) = g_{TJ\_Th}
(1) = \frac{2}{3}.
\end{equation}

By the application of (\ref{eq21}) with (\ref{eq102}) we get
(\ref{eq100}).
\end{proof}

\begin{proposition} \label{pro48} We have
\begin{equation}
\label{eq103}
D_{Th} (P\vert \vert Q) \leqslant 3D_{Jh} (P\vert \vert Q).
\end{equation}
\end{proposition}

\begin{proof} Let us consider
\begin{align}
g_{Th\_Jh} (x) = \frac{{f}''_{Th} (x)}{{f}''_{Jh} (x)} & =
\frac{2\left[ {x^2 + 1 - \sqrt x \left( {x + 1} \right)} \right]}{(x
+ 1)\left( {\sqrt x - 1} \right)^2}, \, x \ne
1\notag\\
& = \frac{2(x + \sqrt x + 1)}{x + 1},\notag
\end{align}

\noindent for all $x \in (0,\infty )$, where ${f}''_{Th} (x)$ and
${f}''_{Jh} (x)$ are as given by (\ref{eq61}) and (\ref{eq64})
respectively.

Calculating the first order derivative of the function $g_{Th\_Jh}
(x)$ with respect to $x$, one gets
\begin{equation}
\label{eq104} {g}'_{Th\_Jh} (x) = - \frac{x - 1}{\sqrt x \left( {x
+ 1} \right)^2} \,
\begin{cases}
 { < 0,} & {x < 1} \\
 { < 0,} & {x > 1} \\
\end{cases}.
\end{equation}

In view of (\ref{eq104}) we conclude that the function $g_{Th\_Th}
(x)$ is increasing in $x \in (0,1)$ and decreasing in $x \in
(1,\infty )$, and hence
\begin{equation}
\label{eq105}
M = \mathop {\sup }\limits_{x \in (0,\infty )} g_{Th\_Jh} (x) = g_{Th\_Jh}
(1) = 3.
\end{equation}

By the application of (\ref{eq21}) with (\ref{eq105}) we get
(\ref{eq103}).
\end{proof}

\begin{remark} \label{not45} In view of Propositions \ref{pro47} and \ref{pro48},
and the inequality (\ref{eq131}) we conclude the following
inequality
\begin{equation}
\label{eq106}
h(P\vert \vert Q) \leqslant \frac{T(P\vert \vert Q) + 2h(P\vert \vert Q)}{3}
\leqslant \frac{1}{8}J(P\vert \vert Q).
\end{equation}
\end{remark}

\begin{proposition} \label{pro49} We have
\begin{equation}
\label{eq107}
D_{Jh} (P\vert \vert Q) \leqslant \frac{1}{12}D_{\Psi \Delta } (P\vert \vert
Q).
\end{equation}
\end{proposition}

\begin{proof} Let us consider
\begin{align}
g_{Jh\_\Psi \Delta } (x) = \frac{{f}''_{Jh} (x)}{{f}''_{\Psi
\Delta } (x)} & = \frac{x\left( {\sqrt x - 1} \right)^2\left( {x +
1} \right)^3}{\left( {x - 1} \right)^2\left( {x^4 + 5x^3 + 12x^2 +
5x + 1} \right)}, \, x \ne 1 \notag\\
& = \frac{x(x + 1)^3}{\left( {\sqrt x + 1} \right)^2\left( {x^4 +
5x^3 + 12x^2 + 5x + 1} \right)}.\notag
\end{align}

\noindent for all $x \in (0,\infty )$, where ${f}''_{Jh} (x)$ and
${f}''_{\Psi \Delta } (x)$ are as given by (\ref{eq64}) and
(\ref{eq59}) respectively.

Calculating the first order derivative of the function
$g_{Jh\_\Psi \Delta } (x)$ with respect to $x$, one gets
\begin{align}
\label{eq108} {g}'_{Jh\_\Psi \Delta } (x) & = - \frac{\left( {\sqrt
x - 1} \right)\left( {x + 1} \right)^2}{\left( {\sqrt x + 1}
\right)^3\left( {x^4 + 5x^3 + 12x^2
+ 5x + 1} \right)^2}\times \\
& \qquad \qquad \times \left[{ {x^5 + 5x^4 + 6x^2(\sqrt x -1)^2 +
5x+1}
} \right. \notag\\
& \qquad \qquad \qquad \left. { + \sqrt x \left( {x^4 + 3x^3 + 4x^2
+ 3x + 1} \right) } \right].\notag
\end{align}

From (\ref{eq108}), one gets
\begin{equation}
\label{eq109} {g}'_{Jh\_\Psi \Delta } (x) \,
\begin{cases}
 { > 0,} & {x < 1} \\
 { < 0,} & {x > 1} \\
\end{cases}.
\end{equation}

In view of (\ref{eq109}) we conclude that the function
$g_{Jh\_\Psi \Delta } (x)$ is increasing in $x \in (0,1)$ and
decreasing in $x \in (1,\infty )$, and hence
\begin{equation}
\label{eq110}
M = \mathop {\sup }\limits_{x \in (0,\infty )} g_{Jh\_\Psi \Delta } (x) =
g_{Jh\_\Psi \Delta } (1) = \frac{1}{12}.
\end{equation}

By the application of (\ref{eq21}) with (\ref{eq110}) we get
(\ref{eq107}).
\end{proof}

\begin{remark} \label{not46} In view of Propositions \ref{pro49},
and the inequality (\ref{eq131}) we conclude the following
inequality
\begin{equation}
\label{eq111} \frac{3}{2}J(P\vert \vert Q) + \frac{1}{4}\Delta
(P\vert \vert Q) \leqslant \frac{1}{16}\Psi (P\vert \vert Q) +
12h(P\vert \vert Q).
\end{equation}
\end{remark}

\begin{proposition} \label{pro410} We have
\begin{equation}
\label{eq112}
D_{\Psi \Delta } (P\vert \vert Q) \leqslant \frac{6}{5}D_{\Psi I} (P\vert
\vert Q).
\end{equation}
\end{proposition}

\begin{proof} Let us consider
\[
g_{\Psi \Delta \_\Psi I} (x) = \frac{{f}''_{\Psi \Delta }
(x)}{{f}''_{\Psi I} (x)} = \frac{x^4 + 5x^3 + 2x^2 + 5x + 1}{(x +
1)^2(x^2 + 3x + 1)}, \, x \in (0,\infty ),
\]

\noindent where ${f}''_{\Psi \Delta } (x)$ and ${f}''_{\Psi I}
(x)$ are as given by (\ref{eq59}) and (\ref{eq58}) respectively.

Calculating the first order derivative of the function $g_{\Psi
\Delta \_\Psi I} (x)$ with respect to $x$, one gets
\begin{equation}
\label{eq113} {g}'_{\Psi \Delta \_\Psi I} (x) = - \frac{4x(x -
1)(2x + 1)(x + 2)}{(x + 1)^3(x^2 + 3x + 1)^2} \,
\begin{cases}
 { > 0,} & {x < 1} \\
 { < 0,} & {x > 1} \\
\end{cases}.
\end{equation}

In view of (\ref{eq113}) we conclude that the function $g_{\Psi
\Delta \_\Psi I} (x)$ is increasing in $x \in (0,1)$ and
decreasing in $x \in (1,\infty )$, and hence
\begin{equation}
\label{eq114}
M = \mathop {\sup }\limits_{x \in (0,\infty )} g_{\Psi \Delta \_\Psi I} (x)
= g_{\Psi \Delta \_\Psi I} (1) = \frac{6}{5}.
\end{equation}

By the application of (\ref{eq21}) with (\ref{eq114}) we get
(\ref{eq112}).
\end{proof}

\begin{remark} \label{not47} In view of Propositions \ref{pro410},
and the inequality (\ref{eq131}) we conclude the following
inequality
\begin{equation}
\label{eq115} I(P\vert \vert Q) \leqslant \frac{1}{6}\left[
{\frac{1}{16}\Psi (P\vert \vert Q) + \frac{5}{4}\Delta (P\vert
\vert Q)} \right] \leqslant \frac{1}{16}\Psi (P\vert \vert Q).
\end{equation}
\end{remark}

\begin{proposition} \label{pro411} We have
\begin{equation}
\label{eq116}
D_{\Psi I} (P\vert \vert Q) \leqslant \frac{10}{9}D_{\Psi h} (P\vert \vert
Q).
\end{equation}
\end{proposition}

\begin{proof} Let us consider
\begin{align}
g_{\Psi I\_\Psi h} (x) = \frac{{f}''_{\Psi I} (x)}{{f}''_{\Psi h}
(x)} & = \frac{(x - 1)^2(x^2 + 3x + 1)}{(x + 1)\left( {x\sqrt x -
1} \right)^2}, \, x \ne 1\notag\\
&  = \frac{\left( {\sqrt x + 1} \right)^2\left( {x^2 + 3x + 1}
\right)}{(x + 1)\left( {x + \sqrt x + 1} \right)^2}.\notag
\end{align}

\noindent for all $x \in (0,\infty )$, where ${f}''_{\Psi I} (x)$
and ${f}''_{\Psi h} (x)$ are as given by (\ref{eq58}) and
(\ref{eq57}) respectively.

Calculating the first order derivative of the function $g_{\Psi
I\_\Psi h} (x)$ with respect to $x$, one gets
\begin{equation}
\label{eq117}
{g}'_{\Psi I\_\Psi h} (x) = - \frac{(x - 1)\left( {3x + \sqrt x + 3}
\right)}{\left( {x + \sqrt x + 1} \right)^3(x + 1)^2}
\begin{cases}
 { > 0,} & {x < 1} \\
 { > 0,} & {x > 1} \\
\end{cases}.
\end{equation}

In view of (\ref{eq117}) we conclude that the function $g_{\Psi
I\_\Psi h} (x)$ is increasing in $x \in (0,1)$ and decreasing in
$x \in (1,\infty )$, and hence
\begin{equation}
\label{eq118}
M = \mathop {\sup }\limits_{x \in (0,\infty )} g_{\Psi I\_\Psi h} (x) =
g_{\Psi I\_\Psi h} (1) = \frac{10}{9}.
\end{equation}

By the application of (\ref{eq21}) with (\ref{eq118}) we get
(\ref{eq116}).
\end{proof}

\begin{remark} \label{not48} In view of Propositions \ref{pro411}, and the inequality (\ref{eq131})
we conclude the following inequality
\begin{equation}
\label{eq119} h(P\vert \vert Q) \leqslant \frac{1}{10}\left[
{\frac{1}{16}\Psi (P\vert \vert Q) + 9I(P\vert \vert Q)} \right]
\leqslant \frac{1}{16}\Psi (P\vert \vert Q).
\end{equation}
\end{remark}

\begin{proposition} \label{pro412} We have
\begin{equation}
\label{eq120}
D_{\Psi h} (P\vert \vert Q) \leqslant \frac{9}{8}D_{\Psi J} (P\vert \vert
Q).
\end{equation}
\end{proposition}

\begin{proof} Let us consider
\[
g_{\Psi h\_\Psi J} (x) = \frac{{f}''_{\Psi h} (x)}{{f}''_{\Psi J}
(x)} = \frac{\left( {x + \sqrt x + 1} \right)^2}{(\sqrt x + 1)^2(x +
1)}, \, x \in (0,\infty ),
\]

\noindent where ${f}''_{\Psi h} (x)$ and ${f}''_{\Psi J} (x)$ are
as given by (\ref{eq57}) and (\ref{eq56}) respectively.

Calculating the first order derivative of the function $g_{\Psi
h\_\Psi J} (x)$ with respect to $x$, one gets
\begin{equation}
\label{eq121} {g}'_{\Psi h\_\Psi J} (x) = - \frac{\left( {\sqrt x -
1} \right)\left( {x + \sqrt x + 1} \right)}{(\sqrt x + 1)^3(x +
1)^2} \,
\begin{cases}
 { > 0,} & {x < 1} \\
 { < 0,} & {x > 1} \\
\end{cases}
\end{equation}

In view of (\ref{eq121}) we conclude that the function $g_{\Psi
h\_\Psi J} (x)$ is monotonically increasing in $x \in (0,1)$ and
decreasing in $x \in (1,\infty )$, and hence
\begin{equation}
\label{eq122} M = \mathop {\sup }\limits_{x \in (0,\infty )} g_{\Psi
h\_\Psi J} (x) = g_{\Psi h\_\Psi J} (1) = \frac{9}{8}.
\end{equation}

By the application of (\ref{eq21}) with (\ref{eq122}) we get
(\ref{eq120}).
\end{proof}

\begin{remark} \label{not49} In view of Propositions \ref{pro412},
and the inequalities (\ref{eq131}) we conclude the following
inequality
\begin{equation}
\label{eq123}
\frac{1}{8}J(P\vert \vert Q) \leqslant \frac{1}{9}\left[ {\frac{1}{16}\Psi
(P\vert \vert Q) + 8h(P\vert \vert Q)} \right] \leqslant \frac{1}{16}\Psi
(P\vert \vert Q).
\end{equation}
\end{remark}

\begin{proposition} \label{pro413} We have
\begin{equation}
\label{eq124}
D_{\Psi J} (P\vert \vert Q) \leqslant \frac{4}{3}D_{\Psi T} (P\vert \vert
Q).
\end{equation}
\end{proposition}

\begin{proof} Let us consider
\[
g_{\Psi J\_\Psi T} (x) = \frac{{f}''_{\Psi J} (x)}{{f}''_{\Psi T}
(x)} = \frac{(x + 1)^2}{x^2 + x + 1}, \, x \in (0,\infty ),
\]

\noindent where ${f}''_{\Psi J} (x)$ and ${f}''_{\Psi T} (x)$ are
as given by (\ref{eq56}) and (\ref{eq55}) respectively.

Calculating the first order derivative of the function $g_{\Psi
J\_\Psi T} (x)$ with respect to $x$, one gets
\begin{equation}
\label{eq125} {g}'_{\Psi J\_\Psi T} (x) = - \frac{(x - 1)(x +
1)}{(x^2 + x + 1)^2} \,
\begin{cases}
 { > 0,} & {x < 1} \\
 { < 0,} & {x > 1} \\
\end{cases}.
\end{equation}

In view of (\ref{eq125}) we conclude that the function $g_{\Psi
J\_\Psi T} (x)$ is increasing in $x \in (0,1)$ and decreasing in
$x \in (1,\infty )$, and hence
\begin{equation}
\label{eq126}
M = \mathop {\sup }\limits_{x \in (0,\infty )} g_{\Psi J\_\Psi T} (x) =
g_{\Psi J\_\Psi T} (1) = \frac{4}{3}.
\end{equation}

By the application of (\ref{eq21}) with (\ref{eq126}) we get
(\ref{eq124}).
\end{proof}

\begin{remark} \label{not410} In view of Propositions \ref{pro413},
and the inequality (\ref{eq131}) we conclude the following
inequality
\begin{equation}
\label{eq127}
T(P\vert \vert Q) \leqslant \frac{1}{32}\left[ {\frac{1}{2}\Psi (P\vert
\vert Q) + 3J(P\vert \vert Q)} \right] \leqslant \frac{1}{16}\Psi (P\vert
\vert Q).
\end{equation}
\end{remark}

Combining (\ref{eq100}), (\ref{eq103}), (\ref{eq107}),
(\ref{eq112}), (\ref{eq116}), (\ref{eq120}) and (\ref{eq124}) we
get (\ref{eq76}). Thus the combination of the Propositions
\ref{pro41}-\ref{pro413} completes the proof of the Theorem
\ref{the41}.

\section{Final Comments}
\begin{itemize}
\item[(i)] In view of inequalities (\ref{eq84}), (\ref{eq88}),
(\ref{eq93}), (\ref{eq100}), (\ref{eq107}) and (\ref{eq127}), we
have the following improvement over the inequality (\ref{eq131}):
\begin{align}
\label{eq128} &\frac{1}{4}\Delta (P\vert \vert Q) \leqslant
I(P\vert \vert Q) \leqslant \frac{2}{3}h(P\vert \vert Q) +
\frac{1}{12}\Delta (P\vert \vert Q) \leqslant h(P\vert \vert Q)\\
 & \qquad \leqslant \frac{1}{16}J(P\vert \vert Q) + \frac{1}{2}I(P\vert \vert Q)
\leqslant \frac{1}{3}T(P\vert \vert Q) + \frac{2}{3}h(P\vert \vert
Q)\notag\\
 & \qquad \qquad \leqslant \frac{1}{8}J(P\vert \vert Q) \leqslant \frac{2}{3}T(P\vert \vert
Q) + \frac{1}{12}\Delta (P\vert \vert Q) \leqslant T(P\vert \vert
Q)\notag\\
 & \qquad \qquad \qquad \leqslant \frac{1}{32}\left[ {\frac{1}{2}\Psi (P\vert \vert Q) + 3J(P\vert
\vert Q)} \right] \leqslant \frac{1}{16}\Psi (P\vert \vert
Q).\notag
\end{align}

\item[(ii)] For simplicity, if we write, the divergence measures
given in (\ref{eq40})-(\ref{eq54}) by $D_1 -D_{15} $ respectively,
then the Theorem \ref{the41} resumes in the following
inequalities:\\
\begin{itemize}
\item[(a)] $D_{15} \leqslant \frac{2}{3}D_{14} \leqslant 2D_{13}
\leqslant D_6 $;

\item[(b)] $D_{15} \leqslant \frac{2}{3}D_{14} \leqslant
\frac{1}{2}D_{12} \leqslant \frac{1}{3}D_9 \leqslant D_6 $;

\item[(c)] $D_6 \leqslant \frac{2}{3}D_7 \leqslant 2D_{10}
\leqslant \frac{1}{6}D_5 \leqslant \frac{1}{5}D_4 \leqslant
\frac{2}{9}D_3 \leqslant \frac{1}{4}D_2 \leqslant \frac{1}{3}D_1
$.\\
\end{itemize}
\item[(iii)] Following the similar lines of the propositions given in section 5,
we can easily prove the following inequality,
\begin{equation}
\label{eq162} D_{\Psi T} (P\vert \vert Q) \leqslant \frac{1}{64}D^
* (P\vert \vert Q).
\end{equation}

\noindent where $D^*(P\vert \vert Q)$ is as given by (\ref{eq154}).

The inequality (\ref{eq162}) together with Theorem \ref{the41} gives
us the following improvement over the inequalities (\ref{eq152}) and
(\ref{eq153}):
\begin{equation}
\label{eq129} D_{J\Delta } (P\vert \vert Q) \leqslant
\frac{1}{2}D_{\Psi J} (P\vert \vert Q) \leqslant \frac{2}{3}D_{\Psi
T} (P\vert \vert Q) \leqslant \frac{1}{96}D^
* (P\vert \vert Q).
\end{equation}

\noindent or equivalently,
\[
D_{12} \leqslant \frac{1}{2}D_2 \leqslant \frac{2}{3}D_1 \leqslant
\frac{1}{96}D^* .
\]

From the inequality (\ref{eq129}) and item (ii)(b)-(c), we observe
that there are many \textit{divergence measures} in between
$D_{J\Delta } (P\vert \vert Q)$ and $D_{\Psi J} (P\vert \vert Q)$.
Thus the inequality (\ref{eq129}) improves the results due to
Dragomir et al. \cite{dsb}.\\

\item[(iv)] The inequalities (\ref{eq111}) and (\ref{eq123}) can
be written as
\begin{align}
\label{eq161} & \frac{1}{8}J(P\vert \vert Q) \leqslant
\frac{1}{12}\left[ {\frac{1}{16}\Psi (P\vert \vert Q) + 12h(P\vert
\vert Q) - \frac{1}{4}\Delta (P\vert \vert Q)} \right]\\
& \qquad \qquad \qquad \leqslant \frac{1}{9}\left[
{\frac{1}{16}\Psi (P\vert \vert Q) + 8h(P\vert \vert Q)} \right]
\leqslant \frac{1}{16} \Psi(P\vert \vert Q).\notag
\end{align}

The middle inequalities of (\ref{eq161}) follow in view of
(\ref{eq75}) and (\ref{eq76}).

\item[(v)] The inequalities (\ref{eq119}) and (\ref{eq123}) can be
written as
\begin{align}
\label{eq130} & h(P\vert \vert Q) \leqslant \frac{1}{10}\left[
{\frac{1}{16}\Psi (P\vert
\vert Q) + 9I(P\vert \vert Q)} \right]\\
 & \qquad \leqslant \frac{1}{9}\left[
{\frac{1}{16}\Psi (P\vert \vert Q) + 8h(P\vert \vert Q)}
\right]\leqslant \frac{1}{16}\Psi (P\vert \vert Q).\notag
\end{align}

\end{itemize}


\begin{thebibliography}{99}
\bibitem{pea} K. PEARSON, On
the criterion that a given system of deviations from the probable in
the case of correlated system of variables is such that it can be
reasonable supposed to have arisen from random sampling,
\textit{Phil. Mag.}, \textbf{50}(1900), 157-172.
\end{thebibliography}
\end{document}